\documentclass[11pt]{amsart} %

\usepackage{soul}
\usepackage{float}

\usepackage{epsfig}

\usepackage{array,amsmath, enumerate,  url, psfrag}
\usepackage{amssymb, amsaddr, fullpage}
\usepackage{graphicx,subfigure}
\usepackage[usenames,dvipsnames]{pstricks}
\usepackage{color}

\restylefloat{figure}
 \makeatother

\newenvironment{Pf}
{\noindent\textbf{Proof.}\ \ }{\hfill $\Box$ \medskip}

\newtheorem{theorem}{Theorem}[section]

\newtheorem{defn}[theorem]{Definition}

\newtheorem{satz}{thm1}

\newtheorem{sat}{thm2}

\newtheorem{thm}[sat]{Theorem}
\newtheorem{lem}[satz]{Lemma}

\title{The strong equitable vertex 1-arboricity of complete bipartite graphs and balanced complete $k$-partite graphs}
\author{
Janejira Laomala \hskip 0.2in
Keaitsuda Maneeruk Nakprasit \hskip 0.2in 
Kittikorn Nakprasit \hskip 0.2in 
Watcharintorn Ruksasakchai}

\begin{document}

\maketitle
\begin{abstract}  
An \emph{equitable $(q, r)$-tree-coloring} of a graph $G$ is a $q$-coloring of $G$ such that the subgraph induced by each color class is a forest of maximum degree at most $r$ and the sizes of any two color classes differ by at most $1.$ 
Let the \emph{strong equitable vertex $r$-arboricity} 
of a graph $G,$ denoted by $va^\equiv_r (G)$, 
be the minimum $p$ such that $G$ has an equitable
$(q, r)$-tree-coloring for every $q\geq p.$ 

The values of $va^\equiv_1 (K_{n,n})$ were investigated by Tao and Lin \cite{Tao} and  Wu, Zhang, and Li \cite{WZL13} where
 exact values of $va^\equiv_1 (K_{n,n})$ were found in some special cases.  
In this paper, we extend their results by giving the exact values  of $va^\equiv_1 (K_{n,n})$ for all cases.  
In the process, we introduce a new function related to 
an equitable coloring and 
obtain a more general result by determining the exact value of each $va^\equiv_1 (K_{m,n})$ and  $va^\equiv_1 (G)$ where   
$G$ is a balanced complete $k$-partite graph $K_{n,\ldots,n}$   
\end{abstract}

\section{Introduction}
Throughout this paper, all graphs considered are finite and simple. 
Let  $V(G)$ and $E(G)$ denote the vertex set and edge set of a graph $G$ respectively.  
Let $K_{n_1,\ldots,n_k}$ be a complete $k$-partite graph 
in which  partite set $X_i$ has size $n_i$ for $1 \leq i \leq k.$ 
Let $K_{k*n}$ denote a complete $k$-partite set with each partite set has size $n,$ 
and naturally we always assume that $k\geq 2.$

An \emph{equitable $k$-coloring} of a graph is a proper vertex
$k$-coloring such that the sizes of every two color classes differ
by at most $1.$  
For a real number $r,$ 
 $\lfloor r \rfloor$  is 
the largest integer number not more than $r,$ 
and $\lceil r \rceil$ is  
the smallest integer number not less than $r.$ 
 A \emph{$k$-set} is a set with $k$ elements. 
A graph $G$ is \emph{equitably $k$-colorable} 
if $G$ admits an equitable $k$-coloring, 
equivalently if $V(G)$ can be partitioned into 
$k$ independent sets where each set is a 
$\lfloor |V(G)|/k \rfloor$-set or 
a $\lceil  |V(G)|/k \rceil$-set.  
The \emph{equitable chromatic number} of a graph $G$ 
is the minimum number $k$ such that $G$ is equitably $k$-colorable. 
In contrast with ordinary proper coloring, 
a graph may have an equitable $k$-coloring  
but has no equitable $(k+1)$-coloring. 
For example, $K_{7,7}$ has an equitable $k$-coloring
for $k=2,4,6$ and $k \ge 8$, but has no equitable $k$-coloring for
$k=3,5,$ and $7$. This leads to the definition of the 
\emph{equitable chromatic threshold} which is 
the minimum $p$ such that $G$ is equitably
$q$-colorable for every $q\geq p,$ 

The topic of equitable coloring, introduced by Meyer \cite{M},  
was motivated by a problem of municipal garbage collection \cite{Garbage}. 
To model the problem, the vertices of the graph are used to represent garbage collection routes. 
If two routes cannot be run in the same day, then two 
corresponding vertices share the same edge. 
Thus the problem of assigning routes for $k$ days satisfying 
the number of routes run on each day are about the same 
can be represented by an equitable $k$-coloring. 
Similar approach can be applied to scheduling \cite{schedule 1, schedule 2, schedule 3, schedule 4}  and 
modeling load balance in parallel memory systems \cite{memory 1, memory 2}.

In \cite{Fan11}, Fan et al introduced an equitable relaxed coloring in which every color class induced a forest with maximum degree at most $1$ and the sizes of any two color classes differ by at most $1.$  
On the basis of the aforementioned research,  Wu, Zhang, and Li \cite{WZL13} introduced 
an \emph{equitable} $(q,r)$\emph{-tree-coloring} of a graph $G$ which is a $q$-coloring of vertices of $G$  such that the subgraph induced by each color class is a forest of maximum degree at most $r$ and the sizes of any two color classes differ by at most one. 
In other words, a graph $G$ has an equitable $(q,r)$-tree-coloring if $V(G)$ 
can be partitioned into $q$ sets such that 
 each set is a 
$\lfloor |V(G)|/q \rfloor$- or  $\lceil  |V(G)|/q \rceil$-set inducing 
a forest with maximum degree at most $r.$
Let the \emph{strong equitable vertex $k$-arboricity}, denoted by $va^\equiv_r (G),$ 
be the minimum $p$ such that $G$ has an equitable
$(q, r)$-tree-coloring for every $q\geq p.$ 

The values of $va^\equiv_1 (K_{n,n})$ were investigated by Tao and Lin \cite{Tao} and  Wu, Zhang, and Li \cite{WZL13} where  
 exact values of $va^\equiv_1 (K_{n,n})$ were found in some special cases.  
In this paper, we extend their results by giving the exact values  of $va^\equiv_1 (K_{n,n})$ for all cases.  
In the process, we introduce a new function related to 
an equitable coloring and 
obtain a furthermore general result by determining the exact value of each $va^\equiv_1 (K_{m,n})$ and  $va^\equiv_1 (G)$ where   
$G$ is a balanced complete $k$-partite graph $K_{n,\ldots,n}$

\section{Helpful Lemmas} 
To find $va^\equiv_1(K_{m,n})$ and $va^\equiv_1(K_{k*n}),$  
we introduce the notion of $p(q: n_1,\ldots, n_k)$ which can be computed 
in  linear time.   

\begin{defn}\label{d1}
Assume that $K_{n_1,\ldots,n_k}$ has an  equitable $q$-coloring. 
Define $p(q: n_1,\ldots, n_k)= \lceil n_1/d \rceil+\cdots+\lceil n_k/d \rceil$  where $d$ is the minimum integer not less than $\lceil (n_1+\cdots +n_k)/q \rceil$  
satisfying at least one of the following conditions:\\ 
\indent (i) there exist $n_i$ and $n_j$ ($i \neq j$) not divisible by $d$;\\  
\indent (ii) there exists $n_i$ with  $n_i/\lfloor n_i/d\rfloor  > d+1.$ 
\end{defn}

\begin{lem}\label{L1} 
Assume $G=K_{n_1,\ldots,n_k}$ has an  equitable $q$-coloring. 
If $p(q: n_1,\ldots, n_k) \leq r \leq q,$  
then $G$ has an equitable $r$-coloring 
\end{lem} 

\begin{Pf} 
Let $p=p(q: n_1, \ldots, n_k)$ and $N=n_1+\cdots +n_k.$ 
We prove by reverse induction that $G$ has an equitable $r$-coloring  when $p \leq r \leq q.$ 
By assumption, $G$ has an equitable  $q$-coloring. 
Assume $G$ has an equitable $r$-coloring where $p < r \leq q.$ 
Let $b= \lceil N/r \rceil.$ 
It follows that, for $1 \leq i \leq k,$ there are nonnegative integers $r_i$ and $s_i$ where $r_i-s_i \geq 0$ such that  
$n_i = (r_i-s_i)b+ s_i(b-1)=r_ib-s_i$ and $r=r_1+\cdots+r_k.$


CASE 1: There exists $j$ with $r_j \neq \lceil n_j/b \rceil.$ 

Note that $n_j = \lceil n_j/b \rceil b - g_j$ for some $g_j$  
satisfying $0 \leq g_j \leq b-1.$ 
From $r_j b -s_j = \lceil n_j/b \rceil b - g_j,$   
we obtain $(r_j -\lceil n_j/b \rceil) b =s_j-g_j.$   
Since $r_j \neq \lceil n_j/b \rceil, 0 \leq g_j \leq b-1,$ 
and  $s_j$ is nonnegative,  
it follows that $s_j -g_j$ is a positive multiple of $b.$ 
Rewrite $n_j = (r_j-s_j)b+s_j(b-1)$ into  $n_j=(r_j-s_j+b-1)b+(s_j-b)(b-1).$ 
Since $s_j -g_j$ is a positive multiple of $b$ and $g_j$ 
is nonnegative, it follows that 
$s_j-b$ is nonnegative. 
Thus we can partition $X_j$ into  $r_j-s_j+b-1$ color classes of size $b$ and 
$s_j-b$ color classes of size $b-1.$ 
That is, we can partition $X_j$ into $r_j-1$ sets of size 
$b$ or $b-1.$
Since we can partition all other $X_i$ into $r_i$ color classes of size $b$ or $b-1$ 
and 
$(\sum_{i \neq j}r_i)+(r_j-1)$ $=(\sum_{i=1}^k r_i)-1=r-1,$  
the graph $G$ has an equitable $(r-1)$-coloring. 

CASE 2: $r_i = \lceil n_i/b \rceil$ for $1 \leq i \leq k.$ 

In this case, $r=\lceil n_1/b \rceil+\cdots+\lceil n_k/b \rceil.$ 
Since $r > p=\lceil n_1/d \rceil+\cdots+\lceil n_k/d \rceil,$   
it follows that  $d > b.$ 
Definition \ref{d1} implies  
 $b$ satisfies neither conditions (i) nor (ii) in the definition. 
Thus we may assume  $n_i=r_i b$ for $2 \leq i \leq k,$  
and  $n_j/\lfloor n_j/b\rfloor  \leq b+1$ for $1 \leq j \leq k.$

SUBCASE 2.1: $n_1 \neq r_1 b.$ 

Then  $b <  n_1/\lfloor n_1/b\rfloor =n_1/(\lceil n_1/b\rceil -1) = n_1/(r_1-1).$  
Since $b$ violates condition (ii), 
we have $n_1/(r_1-1)=n_1/\lfloor n_1/b\rfloor \leq b+1.$ 
Thus  $b<n_1/(r_1-1) \leq b+1.$ 
Consequently, we can partition $n_1$ into $r_1-1$ color classes of size $b$ or $b+1.$ 
Combining with $r_i$ color classes of $X_i$ of size $b$ for $i \geq 2,$ 
we have an equitable $(r-1)$-coloring.  

SUBCASE 2.2 $n_1 = r_1 b.$ 

It follows that $n_i= r_ib$ for  $1\leq i \leq k.$ 
If there is $j$ such that 
$n_j/(r_j-1) \leq b+1,$ then we have an equitable $(r-1)$-coloring as in subcase 2.1. 
Thus we may assume that $n_i/(r_i-1) > b+1$ for $1\leq i \leq k.$ 

We claim that $b+1 \geq d$ and $r_i= \lceil n_i/b\rceil = \lceil n_i/(b+1)\rceil.$   
If the claim holds, we have $r =  \sum_{i=1}^k \lceil n/b \rceil\leq$ 
$\sum_{i=1}^k \lceil n/d \rceil = p$ which contradicts the fact that $r>p.$ 
Thus this situation is impossible.  
To prove the claim, first we suppose to the contrary that there is $n_i$ which is divisible by $b+1.$   
Since $n_i = r_ib,$ we have $r_i =t_i(b+1)$ for some positive integer $t_i.$ 
Consequenly, $n_i/(r_i-1)= t_i(b+1)b/(t_i(b+1)-1) = b+ b/(t_i(b+1)-1) \leq b+1$ 
which contradicts the fact that $n_i/(r_i-1) > b+1.$ 
Thus $n_i$ is not divisible by $b+1$ where  $1\leq i \leq k.$ 
By condition (i) in Definition \ref{d1}, we have $b+1 \geq d.$  
Since $n_i=r_ib$ and $n_i/(r_i-1) > b+1$ where  $1\leq i \leq k,$ 
it follows that $r_i = n_i/b > n_i/(b+1)> r_i-1.$ 
This leads to $r_i=\lceil n_i/b\rceil = \lceil n_i/(b+1)\rceil.$ 
Thus, we have the claim and this completes the proof.       
\end{Pf}

\begin{lem}\label{L2} 
If $G=K_{n_1,\ldots,n_k}$ has an  equitable $q$-coloring, 
then $G$ has no equitable $(p-1)$-coloring 
where $p=p(q: n_1, \ldots, n_k).$ 
\end{lem} 
\begin{Pf}
Recall that $p= \lceil n_1/d \rceil+\cdots+\lceil n_k/d \rceil$ 
where $d$ is as in Definition \ref{d1}. 
Suppose to the contrary that $G$ has an equitable $(p-1)$-coloring.  
Then there is a partite set, say $X_1$ of size $n_1,$  partitioned 
into at most $\lceil n_1/d \rceil -1$ color classes.  
Consequently, there is a color class containing vertices in $X_1$ 
with size at least $d+1.$  
If there is  $X_j$   partitioned 
into at least $\lceil n_j/d \rceil+1$ color classes, 
then there is a color class containing vertices in $X_j$ 
with size at most $d-1.$  So we have two color classes 
with sizes differed by at least $2,$ a contradiction. 
Thus we assume that $X_1$ is partitioned into $\lceil n_1/d \rceil -1$ color classes, 
and $X_j$ is partitioned into $\lceil n_j/d \rceil$ color classes for each $j \in \{2,\ldots,k\}.$ We consider two cases.  

CASE 1: There exists $j \in \{2,\ldots,k\}$ such that  
$n_j$ is not divisible by $d.$  

It follows that  there is a color class containing vertices in $X_j$ 
with size at most $d-1.$  So we have two color classes 
with sizes differed by at least $2,$ a contradiction. 

CASE 2: $X_j$ is partitioned into  exactly $\lceil n_j/d \rceil$ color classes 
and $n_j$ is divisible by $d$ for each $j \in \{2,\ldots,k\}.$ 

Recall that $d$ is as in Definition \ref{d1}. 
It follows that $n_1$ has $n_1/\lfloor n_1/d\rfloor  > d+1.$ 
Since $X_1$ is partitioned into  $\lceil n_1/d \rceil -1$ 
color classes,   there is a  color class containing vertices in $X_1$ with size at least $d+2.$  
But each color class containing vertices in $X_j$ $(j \in \{2,\ldots,k\})$ has 
size $d,$ we have a contradiction. 
Thus $G$ has no equitable $(p-1)$-coloring. 
\end{Pf}

Lemmas \ref{L1} and \ref{L2} yield the following theorem. 
Note that $K_{n_1,\ldots,n_k}$ has an  equitable $r$-coloring 
when $r\geq n_1+\cdots+n_k.$

\begin{thm}\label{T1}
If $K_{n_1,\ldots,n_k}$ has an  equitable $q$-coloring, 
then $p(q: n_1,\ldots, n_k)$ is the minimum $p$ such that 
$G$ has an equitable $r$-coloring for each $r$ satisfying $p \leq r \leq q.$ 
In particular, the equitable chromatic threshold of $K_{n_1,\ldots,n_k}$ is equal to $p(n_1+\cdots+n_k: n_1,\ldots, n_k).$ 
\end{thm}


We use $p(q: n_1,\ldots,n_k)$ to find 
$va^\equiv_1(K_{n_1,\ldots,n_k})$ as follows. 
\begin{lem}\label{L3} 
Let $G=K_{n_1,\ldots,n_k}$ and $N=n_1+\cdots+n_k.$ 
If $G$ has an equitable $q$-coloring where $N/(q-1) \geq 3$ and 
$G$ has an equitable $(r,1)$-tree-coloring for each $r \geq q,$  
then $va^\equiv_1(G) = p(q: n_1,\ldots,n_k).$ 
\end{lem}
\begin{Pf} 
Let $p =p(q: n_1,\ldots,n_k).$ 
From the assumption for $q$ and Definition \ref{d1},   
the graph $G$ has an equitable $(r,1)$-tree coloring for each $r \geq p.$ 
It remains to show that $G$ has no equitable $(p-1,1)$-tree-coloring. 
Suppose to the contrary that $G$ has an equitable $(p-1,1)$-tree-coloring of $G.$ 
Since $p-1 \leq  q-1,$ 
there is a color class of size at least $N/(p-1) \geq N/(q-1) \geq 3$ 
that is not independent.  
Observe that the graph induced by this color class has maximum degree greater than 1 which is a contradiction. 
This completes the proof. 
\end{Pf}

\begin{lem}\label{L4}
If $m=3b+g$ and $n=3c+h$ 
where  $b$ and $c$ are nonnegative integers 
and  $0\leq g, h \leq 2,$   
then $va^\equiv_1(K_{m,n}) \leq b+c+2.$ 
\end{lem}
\begin{Pf} 
Let $N=m+n.$ 
First, consider $q \geq N/2.$ 
We partition $V(K_{m,n})$ into $q$ sets equitably.  
Each resulting set has size not greater than $2,$ so it cannot 
induce a graph with maximum degree more than $1.$    
Thus we have an equitable $(q,1)$-tree-coloring. 

Now assume $b+c+2 \leq  q < N/2.$ 
It follows that $2 < N/q <3.$ 
Thus there are positive integers $r$ and $s$ such that $N = 3r+2s$ and $q= r+s.$  
Since $q \geq b+c+2,$ we have $r \leq b+c.$ 
 
Let $X'_1$  be a $(3b)$-subset of $X_1$ 
and let $X'_2$  be a $(3c)$-subset of $X_2.$ (One of these sets may be empty.) 
Partition $X'_1$ into $b$ $3$-sets 
and partition $X'_2$  into $c$ $3$-sets. 
Since $r \leq b+c,$ we choose $r$ $3$-sets from these $b+c$ sets to 
initiate a new partition.  
Next, we partition set of the remaining $N -3r$ vertices into $s$ $2$-sets.  
Since each of the resulting $3$-sets contains vertices from the same partite set, 
it is an independent set. 
Moreover, each of remaining $2$-sets cannot induce 
a graph with maximum degree more than $1.$ 
Thus this partition is an equitable $(q,1)$-tree-coloring for each $q \geq b+c+2.$ 
\end{Pf}

\section{$va^\equiv_1(K_{m,n})$} 

Let two partite sets of $K_{m,n}$ be $X_1$ and $X_2$ 
where $X_1 = \{u_1,\ldots, u_m\}$ and $X_2=\{v_1, \ldots, v_n\}.$

\begin{lem}\label{L1n} 
$va^\equiv_1(K_{1,n}) = \lceil (n+2)/3\rceil.$ 
\end{lem} 
\begin{Pf} 
Let $r= \lceil (n+2)/3\rceil,$ and   
let $n = 3b+h$ for some integer $h$ where $0 \leq h \leq 2.$ 
Observe that $r = b+1$  when $h =0$ or $1,$ 
and  $r = b+2$  when $h =2.$ 
First, we show that $va^\equiv_1(K_{1,n}) \leq r.$ 
From Lemma \ref{L4}, we know that $K_{1,n}$ has 
an equitable $(q,1)$-tree-coloring for each $q \geq b+2.$  
It remains to show that $K_{1,n}$ has an equitable 
$(b+1,1)$-tree-coloring if $h=0$ or $1.$

From $2 \leq (n+1)/(b+1) = (3b+h+1)/(b+1)<3,$ 
there are a positive integer $r$ and a nonnegative integer $s$ such that 
$n+1=3r+2s$ and $r+s=b+1.$  
Choosing $u_1$ and another vertex in $X_2$ to be in one set,  
and partitioning the set of remaining vertices into $b$ sets equitably, 
we obtain an equitable $(b+1,1)$-tree-coloring. 

Next, we show that $K_{1,n}$ has no equitable $(r-1,1)$-tree-coloring 
to complete the proof. 
Suppose to the contrary that $K_{1,n}$ has an equitable $(r-1,1)$-tree-coloring. 
But  $(n+1)/(r-1) \geq 3.$ 
Consequently, every resulting color class has size at least 3. 
Then a color class containing $u_1$ induces a graph with maximum degree greater than $1,$ a contradiction. This completes the proof. 
 \end{Pf} 

\begin{lem}\label{L2n} 
 $va^\equiv_1(K_{2,n}) = \lceil (n+3)/3\rceil.$ 
\end{lem} 
\begin{Pf} 
Let $r= \lceil (n+3)/3\rceil,$ and  let $n = 3b+h$ for some integer $h$ with $0 \leq h \leq 2.$ 
Observe that $r = b+1$  when $h =0,$  
and  $r = b+2$  when $h =1,$ or $2.$  
First, we show that $va^\equiv_1(K_{2,n}) \leq r.$ 
From Lemma \ref{L4}, we know that $K_{2,n}$ has  
an equitable $(q,1)$-tree-coloring for each $q \geq b+2.$ 
It remains to show that $K_{2,n}$ has an equitable 
$(b+1,1)$-tree-coloring if $n=3b.$

Since $2 \leq  (n+2)/(b+1) = (3b+2)/(b+1)<3,$ 
there are a positive integer $r$ and a nonnegative integer $s$ such that 
$n+2=3r+2s$ and $r+s=b+1.$  
Choosing $u_1$ and $u_2$ to be in one set,  
and partitioning the set of remaining vertices into $b$ sets equitably, 
we obtain an equitable $(b+1,1)$-tree-coloring. 

Next, we show that $K_{2,n}$ has no equitable $(r-1,1)$-tree-coloring 
to complete the proof. 
Suppose to the contrary that $K_{2,n}$ has an equitable $(r-1,1)$-tree-coloring. 
Since  $(2+n)/(r-1) \geq 3,$ every color class has size at least $3.$ 
Thus a color class containing $u_1$ induces a graph with maximum degree greater than 1, 
a contradiction. 
 \end{Pf} 

\begin{lem}\label{L200} 
If $m =3b$ and $n=3c$ for some positive integers $b$ and $c,$  
then   $va^\equiv_1(K_{m,n}) = p(b+c: m,n).$ 
\end{lem} 
\begin{Pf} 
From Lemma \ref{L4}, 
the graph $K_{m,n}$ has an equitable $(q,1)$-tree-coloring for each 
integer $q$ such that $q \geq b+c+2.$ 

Next, we show that $K_{m,n}$ has an equitable $(b+c+1,1)$-tree-coloring. 
We initiate a  partition by assigning three $2$-sets 
$\{u_1,v_1\}, \{u_2,v_2\}, \{u_3,v_3\}$ and then 
partition $X_1 - \{u_1,u_2,u_3\}$ into $b-1$ $3$-sets 
and partition $X_2 - \{v_1,v_2,v_3\}$ into $c-1$ $3$-sets.  
This partition is equivalent to an equitable $(b+c+1)$-tree-coloring as required.  
 
Finally,  we obtain an equitable $(b+c)$-coloring of $K_{m,n}$ 
by partitioning $X_1$ into $b$ $3$-sets and 
partition $X_2$ into $c$ $3$-sets.   
Using Lemma \ref{L3}, we have  $va^\equiv_1(K_{m,n}) =p(b+c:m,n).$ 
\end{Pf} 

\begin{lem}\label{L201} 
If $m =3b$ and $n=3c+1$ for some positive integers $b$ and $c,$  
then   $va^\equiv_1(K_{m,n}) = p(b+c:m,n).$ 
\end{lem} 
\begin{Pf} 
From Lemma \ref{L4}, 
the graph $K_{m,n}$  has  an equitable $(q,1)$-tree-coloring for each 
integer $q$ such that $q \geq b+c+2.$ 

Next, we show that $K_{m,n}$ has an equitable $(b+c+1)$-coloring. 
We initiate a partition by assigning sets $\{v_1,v_2\}, \{v_3,v_4\}$  
and then we  partition $X_1$ into $b$ $3$-sets, 
and partition $X_2 - \{v_1,v_2,v_3,v_4\}$ into $c-1$ $3$-sets 
to obtain an equitable $(b+c+1)$-coloring. 
 
Finally, we  obtain an equitable  $(b+c)$-coloring of $K_{m,n}$ by partitioning 
$X_1$ into $b$ $3$-sets, and partition $X_2$ into  
$c-1$ $3$-sets and one $4$-set.
Using Lemma \ref{L3}, we have  $va^\equiv_1(K_{m,n})= p(b+c:m,n).$ 
\end{Pf}

\begin{lem}\label{L202} 
If $m =3b$ and $n=3c+2$ for some positive integers $b$ and $c,$ 
then   $va^\equiv_1(K_{m,n}) = p(b+c+1:m,n).$ 
\end{lem} 
\begin{Pf} 
From Lemma \ref{L4}, the graph $K_{m,n}$ 
has  an equitable $(q,1)$-tree-coloring for each integer $q$ such that $q \geq b+c+2.$ 
Moreover, we can obtain an equitable  $(b+c+1)$-coloring of $K_{m,n}$   
by partitioning $X_1$ into $b$ $3$-sets, and partitioning $X_2$ into  $c$ 
$3$-sets and one $2$-set . 
Using Lemma \ref{L3}, we have  $va^\equiv_1(K_{m,n})= p(b+c+1:m,n).$ 
\end{Pf}

\begin{lem}\label{L211} 
If $m =3b+1$ and $n=3c+1$ for some positive integers $b$ and $c,$ 
then   $va^\equiv_1(K_{m,n}) = p(b+c:m,n).$ 
\end{lem} 
\begin{Pf} 
From Lemma \ref{L4}, the graph $K_{m,n}$ 
has  an equitable $(q,1)$-tree-coloring for each integer $q$ such that $q \geq  b+c+2.$ 

Next, we show that $K_{m,n}$ has an equitable $(b+c+1)$-tree-coloring. 
We initiate a partition by assigning a set $\{u_1,v_1\}$  
and then we  partition $X_1-\{u_1\}$ into $b$ $3$-sets, 
and partition $X_2 - \{v_1\}$ into $c$ $3$-sets 
to obtain an equitable $(b+c+1)$-coloring. 

To obtain an equitable $(b+c)$-coloring of $K_{m,n},$      
we  partition $X_1$ into $b-1$ $3$-sets and one $4$-set, 
and  we  partition $X_2$ into $c-1$ $3$-sets and one $4$-set. 
Using Lemma \ref{L3}, we have  $va^\equiv_1(K_{m,n})= p(b+c:m,n).$ 
\end{Pf}

\begin{lem}\label{L212} 
If $m =3b+1$ and $n=3c+2$ for some positive integers $b$ and $c,$  
then   $va^\equiv_1(K_{m,n}) = b+c+2.$ 
\end{lem} 
\begin{Pf} 
From Lemma \ref{L4}, the graph $K_{m,n}$  has  an equitable $(q,1)$-tree-coloring for each integer $q$ such that $q \geq b+c+2.$ 

It remains to show that $K_{m,n}$ has no equitable $(b+c+1)$-tree-coloring. 
Suppose to the contrary that $K_{m,n}$ has such coloring. 
Since $(m+n)/(b+c+1)=3,$ each resulting color class has size $3.$  
It follows that each color class is an independent set. 
But $X_1$ can be partitioned into at most $b$ $3$-sets and 
$X_2$ can be partitioned into at most $c$ $3$-sets, a contradiction.  
\end{Pf}

\begin{lem}\label{L222} 
If $m =3b+2$ and $n=3c+2$ for some positive integers $b$ and $c,$ 
then   $va^\equiv_1(K_{m,n}) = b+c+2.$ 
\end{lem} 
\begin{Pf} 
The proof is similar to one of Lemma \ref{L212}. 

\end{Pf} 

From Lemmas \ref{L1n}, \ref{L2n},  
\ref{L200}, \ref{L201}, \ref{L202}, \ref{L211}, \ref{L212}, and \ref{L222}.  
we obtain the following theorem. 
\begin{thm} 
Each $va^\equiv_1(K_{m,n})$ is as in Table 1. 
\end{thm}

\begin{table}[h!]
  \begin{center}
  
    \label{tab:table1}
    \begin{tabular}{l|l|l} 
      $m$ & $n$ & $va_1^\equiv(K_{m,n})$\\
      \hline
      $1$ & $n$ & $\lceil (n+2)/3\rceil$ \\
      $2$ & $n$ & $\lceil(n+3)/3\rceil$\\
      $3b$ & $3c$ & $p(b+c: m,n)$\\
      $3b$ & $3c+1$ & $p(b+c: m,n)$\\
      $3b$ & $3c+2$ & $p(b+c+1: m,n)$\\
      $3b+1$ & $3c+1$ & $p(b+c: m,n)$\\
      $3b+1$ & $3c+2$ & $b+c+2$\\
      $3b+2$ & $3c+2$ & $b+c+2$
      
    \end{tabular}
      \caption{$va_1^\equiv(K_{m,n})$ in terms of $m$ and $n$ where $b$ and $c$ 
    are positive integers}
  \end{center}
\end{table}

\section{$va^\equiv_1(K_{k*n})$} 
For $K_{k*n}$, we denote $p(q: n,\ldots,n)$  by $p(q: k*n).$  

\begin{lem}\label{Lx1} 
If $G = K_{k*1}$ or $K_{k*2}$ and $N = |V(G)|,$  
then $va^\equiv_1(G) = \lceil N/2 \rceil.$   
\end{lem} 
\begin{Pf}
For $q \geq \lceil N/2 \rceil,$ we have $N/q \leq 2.$ 
We obtain an equitable $(q,1)$-tree-coloring by partitioning $V(G)$ into $q$ sets 
equitably. Since each set has size at most $2$, we obtain a desired coloring. 

It remains to  show that $G$ has no equitable $(\lceil N/2 \rceil -1,1)$-tree-coloring. Suppose to the contrary that such coloring exists. 
Then there is a resulting color class with size at least $3$ 
which induces a graph with maximum degree at least $2,$ a contradiction. 
This completes the proof. 
\end{Pf}

\begin{lem}\label{L30} 
If $n =3b$ for some positive integer $b,$   
then   $va^\equiv_1(K_{k*n}) = p(kb: k*n).$ 
\end{lem} 
\begin{Pf} 
For $q \geq 3kb/2,$ we partition $V(K_{k*n})$ into $q$ sets equitably. 
Since each set has size at most $2,$ this partition leads to  
an equitable $(q,1)$-tree-coloring. 

For $kb +1 \leq q \leq 3kb/2 -1,$ then $ 2 < 3kb/q < 3.$  
It follows that there are positive integers $r$ and $s$ such that $3kb = 3r+2s$ 
and $r+s=q.$
Observe that $r < kb.$  
Partition $X_i$ for each $i$ into $b$ $3$-sets to obtain $kb$ $3$-sets.  
Since $r \leq kb,$ we initiate a new partition by choosing $r$ sets from these $kb$ $3$-sets. 
Next, we partition the set of remaining vertices into $s$ sets equitably. 
Since each of $3$-sets contains vertices from the same partite set, 
it is an independent set. 
Moreover, each of remaining sets has size $2.$  
Thus we obtain an equitable $(q,1)$-tree-coloring for each $q \geq kb+1.$ 

One can easily see that $G$ has an equitable $kb$-coloring. 
Using Lemma \ref{L3}, we have  $va^\equiv_1(K_{k*n}) = p(kb: k*n).$ 
\end{Pf}


\begin{lem}\label{L31} 
If $n =3b+1$ for some positive integer $b,$   
then  $va^\equiv_1(K_{2*n}) = p(2b:n,n)$  
and $va^\equiv_1(K_{k*n}) = kb+\lceil k/2 \rceil$ for $k\geq 3.$  
\end{lem} 
\begin{Pf} 
If $k =2,$ the result comes from Lemma \ref{L211}. 
Now assume that $k \geq 3.$ 
For $q \geq (3kb+k)/2,$  we partition $V(K_{k*n})$ into $q$ sets equitably. 
Since each set has size at most $2,$ this partition leads to 
an equitable $(q,1)$-tree-coloring. 

For $kb+\lceil k/2 \rceil \leq q \leq (3kb+k)/2 -1,$ we have $2 < (3kb+k)/q <3.$ 
Thus there are positive integers $r$ and $s$ such that $3kb+k = 3r+2s$ 
and $r+s=q.$ 
Consequently, $3kb+k=3(r+s)-s=3q-s.$  
Since $q \geq kb+\lceil k/2 \rceil,$ 
it follows that $s \geq 3\lceil k/2 \rceil - k.$  
From $3kb+k=3r+2s,$ we have $r \leq kb.$
Choose $X'_i$ to be a $(3b)$-subset of $X_i$ for each $i\in \{1, \ldots, k\}$.  
Partition $X'_i$ for each $i$ into $b$ $3$-sets of size 3. Thus we have $kb$ $3$-sets.  
Since $r \leq kb,$ we initiate a new partition by choosing $r$ $3$-sets from these $kb$ $3$-sets. Next, we partition set of remaining vertices into $s$ sets equitably. 
Since each of  $3$-sets contains vertices from the same partite set, 
it is an independent set. 
Moreover, each of remaining sets has size $2.$ 
Thus we obtain an equitable $(q,1)$-tree-coloring for 
each $q \geq kb+\lceil k/2 \rceil .$

It remains to show that $G$ has no equitable $(kb+\lceil k/2 \rceil-1,1)$-tree-coloring. 
Suppose to the contrary that $G$ has such coloring. 
Consider $k=4.$ This yields the graph $K_{4*n}$ has an equitable $(4b+1,1)$-tree-coloring. 
But $|V(G)|=12b+4=4+(4b)3.$ 
The color classes are $4b$ $3$-sets and one $4$-set.  
Note that each of these color classes must be an independent set. 
But there are at most  $b$ color classes of size $3$ to contain only vertices in $X_i$ where $i\in \{1,2,3, 4\}.$
Thus there are at most $4b$ color classes of size $3,$  
a contradiction. 

Consider $k = 3$ or  $k \geq 5.$  
Since $3kb+k=3(kb-2\lceil k/2 \rceil +k+2)+2(3\lceil k/2 \rceil -k-3),$ 
an equitable $(kb+\lceil k/2 \rceil-1,1)$-tree-coloring yields  
$kb-2\lceil k/2 \rceil +k+2$ color classes of size $3$ and 
$3\lceil k/2 \rceil -k-3)$ color classes of size $2.$ 
Note that a color class of size $3$ must be independent.  
But there are at most  $b$ color classes of size $3$ to contain only vertices in $X_i$ where $i\in \{1,\ldots, k\}.$ 
Thus there are at most $kb$ color classes of size $3,$  
a contradiction. 
\end{Pf} 

\begin{lem}\label{L32} 
Let $n =3b+2$ for some positive integer $b.$  
Then   $va^\equiv_1(K_{k*n}) = kb+k.$ 
\end{lem} 
\begin{Pf} 
If $k=2,$ then we obtain the result  from Lemma \ref{L222}. 
Now assume $k\geq 3.$ 
For $q \geq (3kb+2k)/2,$ we  partition $V(G)$ into $q$ sets equitably. 
Since each set has size at most $2,$ the resulting partition leads to  
an equitable $(q,1)$-tree-coloring. 

For $kb+k \leq q \leq (3kb+2k)/2 -1,$ we have $2 < (3kb+2k)/q < 3.$  
Thus there are positive integers $r$ and $s$ such that 
$3kb+2k = 3r+2s$ and $r+s=q.$ 
Consequently, $3kb+2k=3(r+s)-s=3q-s.$ 
Since $q \geq kb+k,$ we have $s \geq k.$ 
From $3kb+2k=3r+2s,$ we have $r \leq kb.$   
Choose $X'_i$ to be a $(3b)$-subset of $X_i$ where $i\in \{1, \ldots, k\}.$
Partition $X'_i$ for each $i$ into $b$ $3$-sets to obtain $kb$ $3$-sets. 
Since $r \leq kb,$ we can initiate a new partition by choosing $r$ sets from these $kb$ sets.  
Next, we partition the set of the remaining vertices into $s$ $2$-sets.  
Since each of  $3$-sets contains vertices from the same partite set, 
it is an independent set. 
Moreover, each of remaining sets has size $2.$ 
Thus we obtain an equitable $(q,1)$-tree-coloring for each $q \geq kb+k.$ 

It remains to show that $G$ has no equitable $(kb+k-1,1)$-tree-coloring. 
Suppose to the contrary that $G$ has such coloring. 
Since $3kb+2k=3(kb+2)+2(k-3),$ 
an equitable $(kb+k-1,1)$-tree-coloring has $kb+2$ color classes of size $3$  
and $k-3$ color classes of size $2.$ 
Note that a color class of size $3$ must be independent.  
But there are at most $b$ color classes of size  $3$ containing 
only vertices in  $X_i$ where $i\in \{1,\ldots, k\}.$
Thus there are at most $kb$ color classes of size  $3,$  a contradiction.  
\end{Pf} 

From Lemmas \ref{Lx1}, \ref{L30},  \ref{L31},  or \ref{L32}, 
we have the following theorem. 

\begin{thm} 
Each $va^\equiv_1(K_{k*n})$ is as in Table 2. 
\end{thm}

\begin{table}[h!]
  \begin{center}
  
    \label{tab:table1}
    \begin{tabular}{l|l|l} 
       $n$ & $va_1^\equiv(K_{n,n})$ & $va_1^\equiv(K_{k*n})$ where $k\geq 3$\\
      \hline
      $1$ & $1$ & $\lceil k/2\rceil$ \\
      $2$ & $2$ & $k$\\
      $3b$ & $p(2b: n,n)$ &$p(kb: k*n)$ \\
      $3b+1$ & $p(2b:n,n)$ & $kb+\lceil k/2 \rceil$ \\
      $3b+2$ & $2b+2$ & $kb+k$
    
    \end{tabular}
      \caption{$va_1^\equiv(K_{k*n})$ in terms of $n$ where $b$ is a positive integer}
  \end{center}
\end{table}

\section*{Acknowledgments}
The *** author was supported by *********

\end{document}